\documentclass[11pt]{amsart}

\usepackage{epsfig, subfig, amsmath, amssymb, latexsym, amsfonts, xcolor, easybmat, bbm}

\newcommand\eop{{{\hfil \ensuremath \Box}}}

\newenvironment{cor}{\subsection{}{\textbf {Corollary.}}\em}{}
\newenvironment{defn}{\subsection{}{\textbf {Definition.}}\em}{\smallskip}
\newenvironment{eg}{\subsection{}{\textbf {Example.}}}{\smallskip}

\newenvironment{lem}{\subsection{}{\textbf {Lemma.}}\em}{\smallskip}

\newenvironment{thm}{\subsection{}{\textbf {Theorem.}}\em}{\smallskip}
\newenvironment{prob}{\subsection{}{\textbf{Problem.}}}{\smallskip}

\newenvironment{pf}{\noindent{\textbf {Proof.}}} {\begin{flushright}\eop \end{flushright}\smallskip}

\newcommand\sgn{\ensuremath{\mathrm{sgn}\, }}

\newcommand\diag{\ensuremath{\mathrm{diag}\, }}

\newcommand\Root{\ensuremath {\mathrm {root}}}

\newcommand\Dim{\ensuremath {\mathrm {dim}}}
\newcounter{asst}
\newcounter{lab}
\newcounter{asstAA}
\newcounter{asstBB}
\newcounter{asstCC}
\newcounter{asstDD}
\newcounter{asstEE}
\newcounter{asstFF}
\newcounter{asstGG}
\newcounter{asstHH}
\newcounter{asstII}
\newcounter{asstJJ}
\newcounter{asstKK}
\newcounter{asstLL}
\newcounter{asstMM}
\newcounter{suppAA}
\begin{document}


\title{On the spectrum of the symmetric tensor products of certain Hilbert-space operators}


\thanks{${}^1$ Research supported in part by National Natural Science Foundation of China (No.: 12471123)}

\thanks{{\ifcase\month\or Jan.\or Feb.\or March\or April\or May\or
June\or
July\or Aug.\or Sept.\or Oct.\or Nov.\or Dec.\fi\space \number\day,
\number\year}}
\author
	[Y.C.~Yang]{{Yuchi Yang}}
\address
	{School of Mathematics\\
	Jilin University\\
	Changchun 130012\\
	P.R. CHINA}
\email{ycyang23@mails.jlu.edu.cn}

\author
	[Y.H.~Zhang]{{Yuanhang~Zhang${}^1$}}
\address
	{School of Mathematics\\
	Jilin University\\
	Changchun 130012\\
	P.R. CHINA}
\email{zhangyuanhang@jlu.edu.cn}

\begin{abstract}
This paper primarily investigates spectral properties of symmetric tensor products of Hilbert-space operators. For a unilateral weighted shift operator
$S_w$, we present an algorithm to compute the point spectrum of its symmetric and antisymmetric tensor products with the adjoint $S_w^*$. Additionally, we analyze the symmetric tensor product of an injective unilateral weighted shift $S_\alpha$ and a diagonal operator $M$ on $l^2$, demonstrating that its point spectrum must be contained in
$\{0\}$.

\end{abstract}



\maketitle
\markboth{\textsc{  }}{\textsc{}}


\section{Introduction} \label{section1}

In the mid-nineteenth century, tensor products emerged in mathematical literature, notably within Riemann's foundational work on differential geometry  \cite{Rie13, Rie16}. Tensors describe many-body quantum systems \cite{N21}, while symmetric tensors form the basis of general relativity \cite{C04}. Moreover, symmetric tensors
play significant roles in multilinear algebra \cite{G78}, probability \cite{NP05} and statistics \cite{M87}. For an overview of symmetric tensors, we refer to the work of Garcia, O'Loughlin, and Yu \cite{GLY24} and references therein.

In Hilbert-space operator theory, Bhatia \cite{Bha84} precisely evaluated the norm of the derivative for symmetric tensor powers of operators acting on finite-dimensional complex Hilbert spaces. More recently, efforts within the physics community have focused on studying self-adjoint extensions of symmetric tensor products of operators
\cite{IMP14, IP15, LWW21}.

Due to limited literature on symmetric tensor products of non-normal operators in infinite-dimensional Hilbert spaces, recently Garcia, O'Loughlin, and Yu conducted a study of symmetric and antisymmetric tensor products of Hilbert space operators, establishing significant foundational results \cite{GLY24}. Their work particularly focuses on norms and spectra for classes of non-normal operators relevant to function-theoretic operator theory.
The starting point of the present paper is the following problem raised in \cite{GLY24}.

\begin{prob}\cite[Problem 7]{GLY24}\label{P:GLY}
Describe the norm and spectrum of $S_\alpha\odot S_\alpha^*$ and $S_\alpha \wedge S_\alpha^*$, in
which $S_\alpha$ is a weighted shift operator. What can be said if more factors are included?
\end{prob}

For a unilateral weighted shift operator $S_w$, we present an algorithm to compute the point spectrum of its symmetric and antisymmetric tensor products with the adjoint $S_w^*$. This partially resolves the aforementioned problem raised by Garcia, O'Loughlin, and Yu.
Additionally, we characterize the symmetric tensor product of an injective unilateral weighted shift
$S_\alpha$ and a diagonal operator $M$ on $l^2$. Specifically, we demonstrate that the point spectrum of this tensor product must be contained in $\{0\}$, and determine when
it equals to $\{0\}$ precisely. We note that Tian, Wang, and Zhu \cite{TWZ25} addressed norm-related aspects of Problem \ref{P:GLY} very recently.

The paper is organized as follows. Section 2 introduces symmetric and antisymmetric tensor product spaces, the domains for symmetric and antisymmetric tensor products of operators; and presents fundamental results about operator-theoretic properties of symmetric tensor products of
bounded operators. In Section 3, we investigate spectral properties of symmetric and antisymmetric tensor products of forward and backward unilateral weighted shift operators. The concluding section analyzes the point spectrum  and norms for symmetric tensor products of unilateral weighted shifts and diagonal operators.

%
%
%



\section{Preliminaries} \label{section2}



\subsection{}
Let $\mathbb{N}$ be the set of all natural numbers, and $\mathbb{N}_0$ denote $\mathbb{N}\cup \{0\}$. Let $\mathbb{C}$ be the complex field, and $\mathbb{D}$ be the open
unit disc in $\mathbb{C}$. Let $\mathcal{H}$ be a complex separable infinite dimensional Hilbert space endowed with the inner product $\left \langle\cdot,\cdot\right \rangle$. Let $\mathcal{B(H)}$ denote the algebra of all bounded linear operators acting on $\mathcal{H}$. A $subspace$ of $\mathcal{H}$ is a linear manifold which is closed in the norm topology. If $\mathcal{L}\subset\mathcal{H}$, then the $span$ of $\mathcal{L}$, denote by $\vee\mathcal{L}$ is the intersection of all subspaces containing $\mathcal{L}$. Let $\Dim\mathcal{L}$ denote the dimension of $\mathcal{L}$. For $T\in \mathcal{B}(\mathcal{H})$, denote the spectrum, the spectral radius, the point spectrum, the kernel of $T$ by $\sigma(T)$, $r(T)$, $\sigma_{p}(T)$, $\ker T$ respectively. Let $\{e_{0}, e_{1}, \cdots\}$ be the standard orthonormal basis of $\ell^{2}$ and $\{\alpha_i\}_{i=0}^\infty$ be a bounded complex sequence, we define the unilateral weighted shift operator $S_\alpha(e_{i})=\alpha_{i}e_{i+1}~\textup{for}~i=0,1,2,\cdots$.
We refer to the readers to \cite{Shi} for a detailed study of unilateral weighted shift operators.
Up to Hilbert space isomorphism, $\ell^{2}$, $\mathcal{H}$, $H^{2}(\mathbb{D})$ are the same (the definition of $H^{2}(\mathbb{D})$ is in \cite[chapter 5]{GMR23}), we could also define the unilateral weighted shift operators in $\mathcal{B}(\mathcal{H})$, $\mathcal{B}(H^{2}(\mathbb{D}))$ similarly. Given a polynomial $p(x)$, let $\Root(p(x))$ denote the set of the roots of $p(x)$. For $a\in\mathbb{R}$, let $\left\lfloor a\right\rfloor$ denote the largest integer smaller than or equal to $a$.

\begin{defn} \label{defn2.02} \cite{GLY24}
For $u_{1},~u_{2},~\cdots,~u_{n}\in\mathcal{H}$, the simple tensor $u_{1}\otimes u_{2}\otimes\dots \otimes u_{n}:\mathcal{H}^{n}\to\mathbb{C} $ acts as follows:$$(u_{1}\otimes u_{2}\otimes\cdots\otimes u_{n})(v_{1},v_{2},\cdots,v_{n})=\left \langle  u_{1},v_{1}\right \rangle \left \langle  u_{2},v_{2}\right \rangle\cdots\left \langle  u_{n},v_{n}\right \rangle.$$
Let $\mathcal{H}^{\hat{\otimes}n}$ denote the $\mathbb{C}$-vector space spanned by the simple tensors. Actually, from \cite{Sim15} we know that there is a unique inner product on $\mathcal{H}^{\hat{\otimes}n}$ such that
$$\left \langle u_{1}\otimes u_{2}\otimes\cdots\otimes u_{n},v_{1}\otimes v_{2}\otimes\cdots\otimes v_{n} \right \rangle =\left \langle  u_{1},v_{1}\right \rangle \left \langle  u_{2},v_{2}\right \rangle\cdots\left \langle  u_{n},v_{n}\right \rangle$$
for all $u_{1},u_{2},\cdots,u_{n},v_{1},v_{2},\cdots,v_{n}\in\mathcal{H}$.
(With the convention $\mathcal{H}^{\hat{\otimes}0}:=\mathbb{C}$.) For $n=1,2,\cdots$, $\mathcal{H}^{\otimes n}$ denotes the completion of $\mathcal{H}^{\hat{\otimes}n}$ with respect to the inner product.
\end{defn}

\begin{defn} \label{defn2.03} \cite{GLY24}
Let $\sum_{n}$ be the group of permutations of $\left\{1,\cdots,n\right\}$. For all $\pi\in\sum_{n}$ and $u_{1},u_{2},\cdots,u_{n}\in\mathcal{H}$, define $$\hat{\pi}(u_{1}\otimes u_{2}\otimes\cdots\otimes u_{n})=u_{\pi(1)}\otimes u_{\pi(2)}\otimes\cdots\otimes u_{\pi(n)}.$$
\end{defn}	

The density of the span of the simple tensors ensures that $\hat{\pi}$ extends to a bounded linear map on $\mathcal{H}^{\otimes n}$. Then, we will define symmetric and antisymmetric tensor products of Hilbert spaces.
\begin{defn} \label{defn2.04} \cite{GLY24}
Let  $\sgn\pi$ denote the sign of a permutation $\pi\in\sum_{n}$.

\textup{(a)} Let $\mathcal{H}^{\odot 1}=\mathcal{H}$ and $\mathcal{H}^{\odot n}=\left\{v\in\mathcal{H}^{\otimes n}:\hat{\pi}(v)=v~for~all~\pi \in\sum_{n}\right\}~for~n\ge 2.$

\textup{(b)} Let $\mathcal{H}^{\wedge 1}=\left\{0\right\}~and~\mathcal{H}^{\wedge n}=\left\{v\in \mathcal{H}^{\otimes n}:\hat{\pi}(v)=(-1)^{\sgn \pi}v~for~all~\pi\in\sum_{n}\right\}~for~n\ge 2.$
\end{defn}
%
\begin{lem} \label{lem2.05}
\cite[Example 2.5]{GLY24} Let $H^{2}(\mathbb{D})$ denote the Hardy space on $\mathbb{D}$, let $H^{2}(\mathbb{D}^{2})$ denote the Hardy space on the bidisk $\mathbb{D}^{2}$. Then $H^{2}(\mathbb{D})\otimes H^{2}(\mathbb{D})~is~isometrically~isomorphic~to~H^{2}(\mathbb{D}^{2})$.
\end{lem}

\begin{pf}
The $1,z,z^{2},\cdots$ are an orthonormal basis of $H^{2}(\mathbb{D})$, so the simple tensor $z^{i}\otimes \omega^{j}~\textup{for}~i,j=0,1,\cdots$ are an orthonormal basis of $H^{2}(\mathbb{D})\otimes H^{2}(\mathbb{D})$. Note that the map $z^{i}\otimes \omega^{j}\mapsto z^{i}\omega^{j}$ is a unitary, and hence $H^{2}(\mathbb{D})\otimes H^{2}(\mathbb{D})$ is isometrically isomorphism to $H^{2}(\mathbb{D}^{2})$.
\end{pf}	

For a overview of $H^{2}(\mathbb{D}^{2})$ and the operator theory on it, we refer the paper of Douglas and Yang \cite{DG00}. In view of Definition ~\ref{defn2.04} and Lemma ~\ref{lem2.05}, we can identify $H^{2}(\mathbb{D})\odot H^{2}(\mathbb{D})~and~H^{2}(\mathbb{D})\wedge H^{2}(\mathbb{D})$ with
$$H^{2}_{sym}(\mathbb{D}^{2}):=\left\{f(z,\omega)\in H^{2}(\mathbb{D}^{2}):f(z,\omega)=f(\omega,z)~for~all~z,~\omega\in\mathbb{D}\right\},$$
$$H^{2}_{asym}(\mathbb{D}^{2}):=\left\{f(z,\omega)\in H^{2}(\mathbb{D}^{2}):f(z,\omega)=-f(\omega,z)~for~all~z,~\omega\in\mathbb{D}\right\}.$$

\begin{defn} \label{defn2.06} \cite{GLY24}
Let $v_{1},~v_{2},~\cdots,~v_{n}\in\mathcal{H}$,
$$v_{1}\odot v_{2}\odot\cdots\odot v_{n}:=S_{n}(v_{1}\otimes v_{2}\otimes\cdots\otimes v_{n}),$$
$$v_{1}\wedge v_{2}\wedge\cdots\wedge v_{n}:=A_{n}(v_{1}\otimes v_{2}\otimes\cdots\otimes v_{n});$$
where $S_{n},~A_{n}$ are defined in \cite[Definition 2.8.]{GLY24}.
$$S_{n}:=\frac{1}{n!}\sum_{\pi\in\sum_{n}}\hat{\pi},~A_{n}:=\frac{1}{n!}\sum_{\pi\in\sum_{n}}\sgn(\pi)\hat{\pi}.$$
In particular, when $n=2$,
$$v_{1}\odot v_{2}:=\frac{1}{2}(v_{1}\otimes v_{2}+v_{2}\otimes v_{1}),~v_{1}\wedge v_{2}:=\frac{1}{2}( v_{1}\otimes v_{2}-v_{2}\otimes v_{1}).$$

\end{defn}
The following lemmas will be used in Section 3.
\begin{lem} \label{lem2.07} \cite{GLY24}
If $\left\{e_0,e_{1},e_{2},e_{3},\cdots\right\}$ is an orthonormal basis of $\mathcal{H}$, then

$\textup{(a)}~\sqrt{2}(e_{i}\odot e_{j})~for~i<j~and~e_{i}\odot e_{i}~for~i\ge 0~form~an~orthonormal~basis~for~\mathcal{H}^{\odot 2},$

$\textup{(b)}~\sqrt{2}(e_{i}\wedge e_{j})~for~i<j~form~an~orthonormal~basis~for~\mathcal{H}^{\wedge 2},$

$\textup{(c)}~\mathcal{H}^{\otimes 2}=\mathcal{H}^{\odot 2}\oplus \mathcal{H}^{\wedge 2}~is ~an ~orthonormal ~decomposition.$
\end{lem}


\begin{lem} \label{lem2.08} \cite{GLY24}
If $\underset{0\leq i\leq j<\infty}\sum|a_{i,j}|^{2}<\infty ,~then~\underset{0\leq i\leq j<\infty}\sum a_{i,j}e_{i}\odot e_{j}\in \mathcal{H}\odot \mathcal{H}.$
\end{lem}

Finally, we introduce symmetric and antisymmetric tensor products of operators.
\begin{defn} \label{defn2.09} \cite{GLY24}
Let $n\geq 2$ be an integer, $B_{1},B_{2},\cdots,B_{n}\in\mathcal{B}(\mathcal{H})$. Then $B_{1}\odot B_{2}\odot \cdots\odot B_{n}~and~B_{1}\wedge B_{2}\wedge \cdots\wedge B_{n}$ are the restrictions of
$$S_{n}(B_{1},B_{2},\cdots,B_{n})=\frac{1}{n!}\sum_{\pi\in\sum_{n}}\left(B_{\pi(1)}\otimes B_{\pi(2)}\otimes \cdots\otimes B_{\pi(n)}\right).$$
to $\mathcal{H}^{\odot n}~and~\mathcal{H}^{\wedge n}$. Relatively, we may write $B^{\odot n}~and~B^{\wedge n}$ instead of $B\odot B\odot \cdots\odot B$ ($n$ times) and $B\wedge B\wedge \cdots\wedge B$ ($n$ times).
\end{defn}

\begin{eg} \label{eg2.10}
\cite[Example 3.5]{GLY24}
In particular, if $B_{1}, B_{2}\in\mathcal{B}(\mathcal{H})$, we have $B_{1}\odot B_{2}$ and $B_{1}\wedge B_{2}$ are the restrictions of $\displaystyle \frac{1}{2}(B_{1}\otimes B_{2}+B_{2}\otimes B_{1})$ to $\mathcal{H}\odot \mathcal{H}$ and $\mathcal{H}\wedge \mathcal{H}$. Then by Lemma \ref{lem2.07},
$$\frac{1}{2}(B_{1}\otimes B_{2}+B_{2}\otimes B_{1})=\begin{pmatrix}
B_{1}\odot B_{2}&0\\
0&B_{1}\wedge B_{2}\end{pmatrix}:\begin{pmatrix}
\mathcal{H}\odot\mathcal{H}\\
\mathcal{H}\wedge\mathcal{H}\end{pmatrix}\to\begin{pmatrix}
\mathcal{H}\odot\mathcal{H}\\
\mathcal{H}\wedge\mathcal{H}\end{pmatrix}.
$$
\end{eg}






\section{The spectrum of a unilateral weighted shift operator and its adjoint } \label{section3}


%
%

In \cite{GLY24}, Garcia, O'Loughlin, and Yu study the spectrum of $S\odot S^{*}$, $S\wedge S^{*}$, where $S$ is the unilateral shift operator. In this section, we generalize their results to unilateral weighted shift operators.


\begin{defn} \label{defn3.01}
Assuming that $w(i)$ is real-value function on $\mathbb{N}$, we define the following polynomials inductively. For $m,n\in\mathbb{N}$, define
$$D^{(n)}_{m}(x)=xD^{(n)}_{m-1}(x)-\frac{1}{4}w^{2}(m-2)w^{2}(2n-m)D^{(n)}_{m-2}(x),~3\leq m\leq n,$$
$$D^{(n)}_{1}(x)=x,~\textup{for}~\textup{all}~ n\ge1,~D^{(n)}_{2}(x)=x^{2}-\frac{1}{4}w^{2}(0)w^{2}(2n-2),~\textup{for}~\textup{all}~ n\ge2;$$
and
$$K^{(n)}_{m}(x)=xK^{(n)}_{m-1}(x)-\frac{1}{4}w^{2}(m-2)w^{2}(2n-m+1)K^{(n)}_{m-2}(x),~3\leq m\leq n,$$
$$K^{(n)}_{1}(x)=x,~\textup{for}~\textup{all}~ n\ge1,~K^{(n)}_{2}(x)=x^{2}-\frac{1}{4}w^{2}(0)w^{2}(2n-1),~\textup{for}~\textup{all}~ n\ge2.$$

From Laplace Expansion Theorem and Induction, we know that $D_{n}^{(n)}(x)$ equals
$$
 \begin{vmatrix}
x&-\frac{w(0)w(2n-2)}{2}\\
-\frac{w(0)w(2n-2)}{2}&x&-\frac{w(1)w(2n-3)}{2}\\
 &-\frac{w(1)w(2n-3)}{2}&x&\ddots\\
 & &-\frac{w(2)w(2n-4)}{2}&\ddots&-\frac{w(n-3)w(n+1)}{2}\\
 & & &\ddots&x&-\frac{w(n-2)w(n)}{2}\\
 & & & &-\frac{w(n-2)w(n)}{2}&x
 \end{vmatrix}.$$
$K_{n}^{(n)}(x)$ equals
$$
 \begin{vmatrix}
x&-\frac{w(0)w(2n-1)}{2}\\
-\frac{w(0)w(2n-1)}{2}&x&-\frac{w(1)w(2n-2)}{2}\\
 &-\frac{w(1)w(2n-2)}{2}&x&\ddots\\
 & &-\frac{w(2)w(2n-3)}{2}&\ddots&-\frac{w(n-2)w(n+1)}{2}\\
 & & &\ddots&x&-\frac{w(n-1)w(n)}{2}\\
 & & & &-\frac{w(n-1)w(n)}{2}&x
 \end{vmatrix}.$$

Consequently, define
$$C^{+}_{2n}(x)=\left(x-\frac{1}{2}w^{2}(n-1)\right)D^{(n)}_{n-1}(x)-\frac{1}{4}w^{2}(n-2)w^{2}(n)D^{(n)}_{n-2}(x),$$
$$C^{-}_{2n}(x)=\left(x+\frac{1}{2}w^{2}(n-1)\right)D^{(n)}_{n-1}(x)-\frac{1}{4}w^{2}(n-2)w^{2}(n)D^{(n)}_{n-2}(x);$$
and
$$G^{+}_{2n+1}(x)=xK^{(n)}_{n}(x)-\frac{1}{2}w^{2}(n-1)w^{2}(n)K^{(n)}_{n-1}(x),$$
$$G^{-}_{2n+1}(x)=K^{(n)}_{n}(x).$$

In fact, $C^{+}_{2n}(x)$ is the characteristic polynomial of
$$\begin{pmatrix}
 0 & \frac{w(0)w(2n-2)}{2} & \\
 \frac{w(0)w(2n-2)}{2} & 0 & \ddots\\
  & \ddots &\ddots&\frac{w(n-3)w(n+1)}{2}\\
  & &\frac{w(n-3)w(n+1)}{2} &0&\frac{w(n-2)w(n)}{2}\\
  & & &\frac{w(n-2)w(n)}{2}&\frac{w(n-1)w(n-1)}{2}
\end{pmatrix};$$
$C^{-}_{2n}(x)$ is the characteristic polynomial of
$$\begin{pmatrix}
 0 & \frac{w(0)w(2n-2)}{2} & \\
 \frac{w(0)w(2n-2)}{2} & 0 & \ddots\\
  & \ddots &\ddots&\frac{w(n-3)w(n+1)}{2}\\
  & &\frac{w(n-3)w(n+1)}{2} &0&\frac{w(n-2)w(n)}{2}\\
  & & &\frac{w(n-2)w(n)}{2}&-\frac{w(n-1)w(n-1)}{2}
\end{pmatrix};$$
$G^{+}_{2n+1}(x)$ is the characteristic polynomial of
$$\begin{pmatrix}
 0 & \frac{w(0)w(2n-1)}{2} & \\
 \frac{w(0)w(2n-1)}{2} & 0 & \ddots\\
  & \ddots &\ddots&\frac{w(n-2)w(n+1)}{2}\\
  & &\frac{w(n-2)w(n+1)}{2} &0&\frac{\sqrt{2}w(n-1)w(n)}{2}\\
  & & &\frac{\sqrt{2}w(n-1)w(n)}{2}&0
\end{pmatrix};$$
$G^{-}_{2n+1}(x)$ is the characteristic polynomial of
$$\begin{pmatrix}
 0 & \frac{w(0)w(2n-1)}{2} & \\
 \frac{w(0)w(2n-1)}{2} & 0 & \ddots\\
  & \ddots &\ddots&\frac{w(n-3)w(n+2)}{2}\\
  & &\frac{w(n-3)w(n+2)}{2} &0&\frac{w(n-2)w(n+1)}{2}\\
  & & &\frac{w(n-2)w(n+1)}{2}&0
\end{pmatrix}.$$

\end{defn}

The following lemma is well known. We omit the proof.

\begin{lem} \label{lem3.02}
If $A\in \mathcal{B}(\mathcal{H})$, $\mathcal{H}_i$ be a reducing subspace of $A$, $i=1,2,\cdots$, and
$\mathcal{H}=\bigoplus_{i=1}^\infty \mathcal{H}_{i}$. Then
\[\sigma_{p}(A)=\bigcup_{i=1}^\infty \sigma_{p}(A|_{\mathcal{H}_{i}}).\]
\end{lem}

From \cite[Chapter 2]{Shi}, we know that a unilateral weighted shift operator with weighted sequence $\{w(i)\}_{i=0}^{\infty}$ is unitarily equivalent to the unilateral weighted shift operator with weighted sequence $\{|w(i)|\}_{i=0}^\infty$, thus we may assume that $w(i)\ge0$ from now on without loss of generality. Below the brackets $\{\{$ and $\}\}$ indicate a $\it{multiset}$;
 that is, a set that permits multiplicity.
Motivated by \cite[Problem 7]{GLY24}, the following theorem is a generalization  of  \cite[Theorem 8.1]{GLY24}.

\begin{thm} \label{thm3.03}
Let $\{w(i)\}_{i=0}^{\infty}$ be a bounded sequence of nonnegative real numbers, $S_w(z^{i})=w(i)z^{i+1}$ be a unilateral weighted shift operator on $H^{2}(\mathbb{D})$.
Let $C^{+}_{2n}(x), C^{-}_{2n}(x)$, $G^{+}_{2n+1}(x)$, $G^{-}_{2n+1}(x)$ be as in Definition \ref{defn3.01},
set
$$\Root\left(C^{+}_{2n}(x)\right) = \left\{p_{i}^{(2n)},~1\leq i\leq n\right\},~\Root\left(C^{-}_{2n}(x)\right) = \left\{q_{i}^{(2n)},~1\leq i\leq n\right\},$$
$$\Root\left(G^{+}_{2n+1}(x)\right) = \left\{p_{i}^{(2n+1)},~1\leq i\leq n+1\right\},~\Root\left(G^{-}_{2n+1}(x)\right) = \left\{q_{i}^{(2n+1)},~1\leq i\leq n\right\},$$ \\counting multiplicity.
Then
$$\sigma_{p}(S_w\odot S_w^{*})=\left\{\left\{p_{i}^{(j)}:j\ge 2~\textup{and}~1\leq i\leq \left\lfloor\frac{j+1}{2}\right\rfloor\right\}\right\}\cup\left\{0\right\},$$
$$\sigma_{p}(S_w\wedge S_w^{*})=\left\{\left\{q_{i}^{(j)}:j\ge 2~\textup{and}~1\leq i\leq\left\lfloor\frac{j}{2}\right\rfloor\right\}\right\}.$$
\end{thm}

\begin{pf} Recall that, $\left\{z^{i}\omega^{j},~i,~j\in\mathbb{N}_0\right\}$ form an orthonormal basis of $H^{2}(\mathbb{D}^{2})$, and $\left\{z^{i},~i\in\mathbb{N}_0\right\}$ form an orthonormal basis of $H^{2}(\mathbb{D})$. Set
$$T=\frac{1}{2}(S_w\otimes S_w^{*}+S_w^{*}\otimes S_w).$$
From Lemma ~\ref{lem2.05}, we can consider $T$ in $\mathcal{B}(H^{2}(\mathbb{D}^{2}))$. Then
$$T(z^{i}\omega^{j})=\begin{cases}\frac{1}{2}w(i)w(j-1)z^{i+1}\omega^{j-1}+\frac{1}{2}w(i-1)w(j)z^{i-1}\omega^{j+1}&\text{if}~~i,~j\ge 1,\\ \frac{1}{2}w(i)w(j-1)z^{i+1}\omega^{j-1}&\text{if}~~i=0,~j\ge 1,\\ \frac{1}{2}w(i-1)w(j)z^{i-1}\omega^{j+1}&\text{if}~~i\ge 1,~j=0,\\ 0&\text{if}~~i=j=0.\end{cases}$$
Define $V_{0}=V_{0}^{+}=\vee\left\{1\right\},~V_{0}^{-}=\left\{0\right\}$. When $k\ge 1,$ set
$$V_{k}=\vee\left\{z^{i}\omega^{k-i}:0\leq i\leq k\right\},~\textup{and}~\textup{hence}~\Dim V_{k}=k+1.$$
Similarly, set
$$V_{k}^{+}=\vee\left\{z^{i}\omega^{k-i}+z^{k-i}\omega^{i}:0\leq i\leq \left\lfloor\frac{k}{2}\right\rfloor\right\},~\textup{and}~\textup{hence}~\Dim V_{k}^{+}=\left\lfloor\frac{k}{2}\right\rfloor+1.$$
$$V_{k}^{-}=\vee\left\{z^{i}\omega^{k-i}-z^{k-i}\omega^{i}:0\leq i\leq \left\lfloor\frac{k-1}{2}\right\rfloor\right\},~\textup{and}~\textup{hence}~\Dim V_{k}^{-}=\left\lfloor\frac{k-1}{2}\right\rfloor+1.$$
Thus
$$V_{k}=V_{k}^{+}\oplus V_{k}^{-},~\Dim V_{k}=\Dim V_{k}^{+}+\Dim V_{k}^{-},~~k\in \mathbb{N}_0.$$
From Lemma ~\ref{lem2.07}, we have $H^{2}(\mathbb{D}^{2})=H_{sym}^{2}(\mathbb{D}^{2})\oplus H_{asym}^{2}(\mathbb{D}^{2}),~V_{k}=V_{k}^{+}\oplus V_{k}^{-},~k\ge 1$. By \cite[(8.2)]{GLY24}, we know
$$H^{2}(\mathbb{D}^{2})=\underset{k=0}{\overset{\infty}{\oplus}}V_{k},~H^{2}_{sym}(\mathbb{D}^{2})=\underset{k=0}{\overset{\infty}{\oplus}}V_{k}^{+},~H^{2}_{asym}(\mathbb{D}^{2})=\underset{k=1}{\overset{\infty}{\oplus}}V_{k}^{-}.$$

Next we will show that $V_{k}^{+},~V_{k}^{-}$ are the invariant subspaces of $T$.

\textbf{Case 1}, assume that $1\leq k$ is odd.
Let $n\in \mathbb{N}$ satisfies $\left\lfloor\frac{k}{2}\right\rfloor=\frac{k-1}{2}=n-1$.

\begin{enumerate}
\item

Consider the basis $\mathcal{E}_k^+:=\left\{\frac{1}{\sqrt{2}}(z^{k-i}\omega^{i}+z^{i}\omega^{k-i})\right\}_{i=0}^{\frac{k-1}{2}}$ of $V_k^+$. For $i\ge1,$
\begin{equation}
\begin{aligned}
T\left(\frac{1}{\sqrt{2}}(z^{k-i}\omega^{i}+z^{i}\omega^{k-i})\right) =&\frac{1}{2}(S_w\otimes S_w^{*}+S_w^{*}\otimes S_w)\left(\frac{1}{\sqrt{2}}(z^{k-i}\omega^{i}+z^{i}\omega^{k-i})\right)\\
 =&\frac{\sqrt{2}}{4}(w(k-i)w(i-1)z^{k-i+1}\omega^{i-1}+w(i)w(k-i-1)z^{i+1}\omega^{k-i-1}\\
  &+w(k-i-1)w(i)z^{k-i-1}\omega^{i+1}+w(i-1)w(k-i)z^{i-1}\omega^{k-i+1})\\
 =&\frac{w(i-1)w(k-i)}{2}\left(\frac{1}{\sqrt{2}}(z^{i-1}\omega^{k-i+1}+z^{k-i+1}\omega^{i-1})\right)\\
  &+\frac{w(i)w(k-i-1)}{2}\left(\frac{1}{\sqrt{2}}(z^{i+1}\omega^{k-i-1}+z^{k-i-1}\omega^{i+1})\right)\in V_k^+;\nonumber
\end{aligned}
\end{equation}
for $i=0$,
$$T\left(\frac{1}{\sqrt{2}}(z^{k}+\omega^{k})\right)=\frac{w(0)w(k-1)}{2}\left(\frac{1}{\sqrt{2}}(z^{1}\omega^{k-1}+z^{k-1}\omega^{1})\right)\in V_k^+.$$

Hence, $V_k^+$ is an invariant subspace of $T$. Set
$B_{k}^{+}=T|_{V_{k}^{+}}$.
With respect to $\mathcal{E}_k^+$, the tridiagonal matrix
$$B_{k}^{+}=
\begin{pmatrix}
0&\frac{w(0)w(k-1)}{2}\\
\frac{w(0)w(k-1)}{2}&0&\frac{w(1)w(k-2)}{2}\\
 &\frac{w(1)w(k-2)}{2}&0&\ddots\\
 & &\frac{w(2)w(k-3)}{2}&\ddots&\frac{w(\frac{k-5}{2})w(\frac{k+3}{2})}{2}\\
 & & &\ddots&0&\frac{w(\frac{k-3}{2})w(\frac{k+1}{2})}{2}\\
 & & & &\frac{w(\frac{k-3}{2})w(\frac{k+1}{2})}{2}&\frac{w(\frac{k-1}{2})w(\frac{k-1}{2})}{2}
 \end{pmatrix},$$
then $xI-B_{k}^{+}$ equals
 $$
 \begin{pmatrix}
x&-\frac{w(0)w(k-1)}{2}\\
-\frac{w(0)w(k-1)}{2}&x&-\frac{w(1)w(k-2)}{2}\\
 &-\frac{w(1)w(k-2)}{2}&x&\ddots\\
 & &-\frac{w(2)w(k-3)}{2}&\ddots&-\frac{w(\frac{k-5}{2})w(\frac{k+3}{2})}{2}\\
 & & &\ddots&x&-\frac{w(\frac{k-3}{2})w(\frac{k+1}{2})}{2}\\
 & & & &-\frac{w(\frac{k-3}{2})w(\frac{k+1}{2})}{2}&x-\frac{w(\frac{k-1}{2})w(\frac{k-1}{2})}{2}
 \end{pmatrix}.$$

Recall that, $D^{(n)}_{n}(x)$ equals
$$
 \begin{vmatrix}
x&-\frac{w(0)w(2n-2)}{2}\\
-\frac{w(0)w(2n-2)}{2}&x&-\frac{w(1)w(2n-3)}{2}\\
 &-\frac{w(1)w(2n-3)}{2}&x&\ddots\\
 & &-\frac{w(2)w(2n-4)}{2}&\ddots&-\frac{w(n-3)w(n+1)}{2}\\
 & & &\ddots&x&-\frac{w(n-2)w(n)}{2}\\
 & & & &-\frac{w(n-2)w(n)}{2}&x
 \end{vmatrix}.$$
Thus,
$$|xI-B_{k}^{+}|=\left(x-\frac{1}{2}w^{2}(n-1)\right)D^{(n)}_{n-1}(x)-\frac{1}{4}w^{2}(n-2)w^{2}(n)D^{(n)}_{n-2}(x).$$
In other words, $|xI-B_{k}^{+}|=C^{+}_{2n}(x)$. Therefore, $\sigma_p(B_k^+)=\Root\left(C^{+}_{2n}(x)\right)$.

\item
Consider the basis $\mathcal{E}_k^-:=\left\{\frac{1}{\sqrt{2}}(z^{i}\omega^{k-i}-z^{k-i}\omega^{i})\right\}_{i=0}^{\frac{k-1}{2}}$ of $V_k^-$. For $i\ge 1$,
\begin{equation}
\begin{aligned}
T\left(\frac{1}{\sqrt{2}}(z^{i}\omega^{k-i}-z^{k-i}\omega^{i})\right) =&\frac{1}{2}(S_w\otimes S_w^{*}+S_w^{*}\otimes S_w)\left(\frac{1}{\sqrt{2}}(z^{i}\omega^{k-i}-z^{k-i}\omega^{i})\right)\\
 =&\frac{\sqrt{2}}{4}(w(i)w(k-i-1)z^{i+1}\omega^{k-i-1}-w(k-i)w(i-1)z^{k-i+1}\omega^{i-1}\\
 &+w(i-1)w(k-i)z^{i-1}\omega^{k-i+1}-w(k-i-1)w(i)z^{k-i-1}\omega^{i+1})\\
 =&\frac{w(i)w(k-i-1)}{2}\left(\frac{1}{\sqrt{2}}(z^{i+1}\omega^{k-i-1}-z^{k-i-1}\omega^{i+1})\right)\\
 &+\frac{w(i-1)w(k-i)}{2}\left(\frac{1}{\sqrt{2}}(z^{i-1}\omega^{k-i+1}-z^{k-i+1}\omega^{i-1})\right)\in V_{k}^{-};\nonumber
\end{aligned}
\end{equation}
for $i=0$,
$$T\left(\frac{1}{\sqrt{2}}(\omega^{k}-z^{k})\right)=\frac{w(0)w(k-1)}{2}\left(\frac{1}{\sqrt{2}}(z^{1}\omega^{k-1}-z^{k-1}\omega^{1})\right)\in V_k^-.$$
Therefore, $V_k^-$ is an invariant subspace of $T$. Set
$B_{k}^{-}=T|_{V_{k}^{-}}$.
So, with respect to $\mathcal{E}_k^-$, we have $xI-B_{k}^{-}$ equals
$$
 \begin{pmatrix}
x&-\frac{w(0)w(k-1)}{2}\\
-\frac{w(0)w(k-1)}{2}&x&-\frac{w(1)w(k-2)}{2}\\
 &-\frac{w(1)w(k-2)}{2}&x&\ddots\\
 & &-\frac{w(2)w(k-3)}{2}&\ddots&-\frac{w(\frac{k-5}{2})w(\frac{k+3}{2})}{2}\\
 & & &\ddots&x&-\frac{w(\frac{k-3}{2})w(\frac{k+1}{2})}{2}\\
 & & & &-\frac{w(\frac{k-3}{2})w(\frac{k+1}{2})}{2}&x+\frac{w(\frac{k-1}{2})w(\frac{k-1}{2})}{2}
 \end{pmatrix},$$
and hence $$|xI-B_{k}^{-}|=\left(x+\frac{1}{2}w^{2}(n-1)\right)D^{(n)}_{n-1}(x)-\frac{1}{4}w^{2}(n-2)w^{2}(n)D^{(n)}_{n-2}(x).$$
In other words, $|xI-B_{k}^{-}|=C^{-}_{2n}(x).$ So $\sigma(B_k^-)=\Root\left(C^{-}_{2n}(x)\right)$.
\end{enumerate}

\textbf{Case 2}, assume that $2\leq k=2n$ is even.

\begin{enumerate}
\item
Consider the basis $\mathcal{E}_k^+:=\left\{\frac{1}{\sqrt{2}}(z^{i}\omega^{k-i}+z^{k-i}\omega^{i})\right\}_{i=0}^{\frac{k}{2}-1}\cup \left\{z^{\frac{k}{2}}\omega^{\frac{k}{2}}\right\}$ of $V_k^+$.
As done in Case 1 (1), a tedious but straightforward calculation shows that $V_k^+$ is an invariant subspace of $T$. Set $B_k^+=T|_{V_k^+}$. Then with respect to $\mathcal{E}_k^+$,
$$B_{k}^{+}=\begin{pmatrix}
0&\frac{w(0)w(2n-1)}{2}\\
\frac{w(0)w(2n-1)}{2}&0&\frac{w(1)w(2n-2)}{2}\\
 &\frac{w(1)w(2n-2)}{2}&0&\ddots\\
 & &\frac{w(2)w(2n-3)}{2}&\ddots&\frac{w(n-2)w(n+1)}{2}\\
 & & &\ddots&0&\frac{\sqrt{2}w(n-1)w(n)}{2}\\
 & & & &\frac{\sqrt{2}w(n-1)w(n)}{2}&0
 \end{pmatrix}.$$
Thus, $|xI-B_{k}^{+}|$ equals
$$
 \begin{vmatrix}
x&-\frac{w(0)w(2n-1)}{2}\\
-\frac{w(0)w(2n-1)}{2}&x&-\frac{w(1)w(2n-2)}{2}\\
 &-\frac{w(1)w(2n-2)}{2}&x&\ddots\\
 & &-\frac{w(2)w(2n-3)}{2}&\ddots&-\frac{w(n-2)w(n+1)}{2}\\
 & & &\ddots&x&-\frac{\sqrt{2}w(n-1)w(n)}{2}\\
 & & & &-\frac{\sqrt{2}w(n-1)w(n)}{2}&x
 \end{vmatrix}.$$
Hence, $|xI-B_{k}^{+}|=G^{+}_{2n+1}(x)$. Consequently, $\sigma_p(B_k^+)=\Root\left(G^{+}_{2n+1}(x)\right)$.
%

\item
Similarly, consider the basis $\mathcal{E}_k^-:=\left\{\frac{1}{\sqrt{2}}(z^{i}\omega^{k-i}-z^{k-i}\omega^{i})\right\}_{i=0}^{\frac{k}{2}-1}$ of $V_k^-$.
A second tedious but direct calculation shows that $V_k^-$ is an invariant subspace of $T$. Set $B_k^-=T|_{V_k^-}$. Then
with respect to $\mathcal{E}_k^-$,
$$B_{k}^{-}=\begin{pmatrix}
0&\frac{w(0)w(2n-1)}{2}\\
\frac{w(0)w(2n-1)}{2}&0&\frac{w(1)w(2n-2)}{2}\\
 &\frac{w(1)w(2n-2)}{2}&0&\ddots\\
 & &\frac{w(2)w(2n-3)}{2}&\ddots&\frac{w(n-3)w(n+2)}{2}\\
 & & &\ddots&0&\frac{w(n-2)w(n+1)}{2}\\
 & & & &\frac{w(n-2)w(n+1)}{2}&0
 \end{pmatrix}.$$
Thus, $|xI-B_{k}^{-}|$ equals
$$
 \begin{vmatrix}
x&-\frac{w(0)w(2n-1)}{2}\\
-\frac{w(0)w(2n-1)}{2}&x&-\frac{w(1)w(2n-2)}{2}\\
 &-\frac{w(1)w(2n-2)}{2}&x&\ddots\\
 & &-\frac{w(2)w(2n-3)}{2}&\ddots&-\frac{w(n-3)w(n+2)}{2}\\
 & & &\ddots&x&-\frac{w(n-2)w(n+1)}{2}\\
 & & & &-\frac{w(n-2)w(n+1)}{2}&x
 \end{vmatrix}.$$
That is, $|xI-B_{k}^{-}|=G^{-}_{2n+1}(x)$. Thus, $\sigma_p(B_k^-)=\Root\left(G^{-}_{2n+1}(x)\right)$.
\end{enumerate}

Since $T$ is self-adjoint, $V_{k}^{+},~V_{k}^{-}$ are the reducing subspaces of $T$. As $T(1)=0\in V_{0}^{+}$, $\sigma_p(T|_{V_0^+})=\{0\}$.
Recall that $H^{2}_{sym}(\mathbb{D}^{2})=\underset{k=0}{\overset{\infty}{\oplus}}V_{k}^{+},~H^{2}_{asym}(\mathbb{D}^{2})=\underset{k=1}{\overset{\infty}{\oplus}}V_{k}^{-}$,
and
$$\Root\left(C^{+}_{2n}(x)\right)=\left\{p_{i}^{(2n)},~1\leq i\leq n\right\},~\Root\left(C^{-}_{2n}(x)\right)=\left\{q_{i}^{(2n)},~1\leq i\leq n\right\},$$
$$\Root\left(G^{+}_{2n+1}(x)\right)=\left\{p_{i}^{(2n+1)},~1\leq i\leq n+1\right\},~\Root\left(G^{-}_{2n+1}(x)\right)=\left\{q_{i}^{(2n+1)},~1\leq i\leq n\right\},$$ \\counting multiplicity.
By Lemma \ref{lem3.02},
$$\sigma_{p}(S_w\odot S_w^{*})=\left\{\left\{p_{i}^{(j)}:j\ge 2~\textup{and}~1\leq i\leq \left\lfloor\frac{j+1}{2}\right\rfloor\right\}\right\}\cup\{0\},$$
and
$$\sigma_{p}(S_w\wedge S_w^{*})=\left\{\left\{q_{i}^{(j)}:j\ge 2~\textup{and}~1\leq i\leq\left\lfloor\frac{j}{2}\right\rfloor\right\}\right\}.$$
\end{pf}

\subsection{}
Why is this not the whole answer? Because how do you give the explicit forms of $C_{2n}^+(x), C_{2n}^-(x)$ and likewise $G_{2n+1}^+(x), G_{2n+1}^-(x)$ and calculate their roots.
The resulting problem is strictly a special function theoretic one.

\bigskip

We next consider the case that $w(i)=\frac{1}{a^{i}},~a\ge1$ in Theorem ~\ref{thm3.03}, as well as we can say more about the spectrum. When $a>1$, $S_w$ is a Donoghue operator (see \cite[Chapter 4, Section 4]{RR73}). The following corollary is an application of Theorem ~\ref{thm3.03}.

%
\begin{cor} \label{cor3.04}
Let $a\geq 1$,  $w(i)=\frac{1}{a^{i}}$, $i\in \mathbb{N}_0$, and $S_w$ denote the unilateral weighted shift acting on $H^2(\mathbb{D})$ with $S_w(z^i)=w(i)z^{i+1}$, $i\in \mathbb{N}_0$. Then
$$\sigma_{p}(S_w\odot S_w^{*})=\left\{\left\{\frac{1}{a^{k-1}}\cos\left(\frac{(2j-1)\pi}{k+2}\right):k\ge 0,~1\leq j\leq \left \lfloor \frac{k+2}{2} \right \rfloor\right\}\right\},$$
$$\sigma_{p}(S_w\wedge S_w^{*})=\left\{\left\{\frac{1}{a^{k-1}}\cos\left(\frac{2j\pi}{k+2}\right):k\ge 1,~1\leq j\leq \left \lfloor \frac{k+1}{2} \right \rfloor\right\}\right\},$$
with the eigenvalues in these multisets repeated by multiplicity. In particular, we have
$$\Dim\ker(S_w\odot S_w^{*})=\infty,~\Dim\ker(S_w\wedge S_w^{*})=\infty.$$
 If $a>1$,
$$\sigma(S_w\odot S_w^{*})=\sigma_{p}(S_w\odot S_w^{*}),~\sigma(S_w\wedge S_w^{*})=\sigma_{p}(S_w\wedge S_w^{*}).$$
\end{cor}

\begin{pf} We first show that when $a>1$, $S_w\odot S_w^{*}$ and $S_w\wedge S_w^{*}$ are compact self-adjoint operators. As
$$T=\frac{1}{2}(S_w\otimes S_w^{*}+S_w^{*}\otimes S_w)=T^{*},$$
$T$ is self-adjoint. Since
$$S_w\odot S_w^{*}=T|_{H^{2}_{sym}(\mathbb{D}^{2})},~S_w\wedge S_w^{*}=T|_{H^{2}_{asym}(\mathbb{D}^{2})},$$
$S_w\odot S_w^{*}$ and $S_w\wedge S_w^{*}$ are self-adjoint. In order to show that $S_w\odot S_w^{*}$ and $S_w\wedge S_w^{*}$ are compact, it suffices to show that $T$ is a compact operator.
As $w(i)$ tend to $0$, $S_w$ and $S_w^{*}$ are compact.
Then there exist finite rank operators $ A_{n},B_{n}$ such that $\|A_{n}-S_w\|\to 0,~\|B_{n}-S_w^{*}\|\to 0.$
Since $A_{n}$ and $B_{n}$ are finite rank operators, it is straightforward to check that $A_{n}\otimes B_{n}$ is also a finite rank operator.
Now
\begin{equation}
\begin{aligned}
\|A_{n}\otimes B_{n}-S_w\otimes S_w^{*}\| &=\|A_{n}\otimes B_{n}-S_w\otimes B_{n}+S_w\otimes B_{n}-S_w\otimes S_w^{*}\|\\
 &\leq \|(A_{n}-S_w)\otimes B_{n}\|+\|S_w\otimes(B_{n}-S_w^{*})\|\\
 &=\|A_{n}-S_w\|\|B_{n}\|+\|S\|\|B_{n}-S_w^{*}\|\to 0.\nonumber
\end{aligned}
\end{equation}
So $S_w\otimes S_w^{*}$ is compact when $a>1$, and hence $T$ is a compact operator.

\textbf{Case 1}, when $k\geq 1$ is odd, set $\displaystyle n=\frac{k+1}{2}$. We know that
$C^{+}_{2n}(x)$ is the characteristic polynomial of
$$\frac{1}{a^{k-1}}\begin{pmatrix}
 0 & \frac{1}{2} & \\
 \frac{1}{2} & 0 & \ddots\\
  & \ddots &\ddots&\frac{1}{2}\\
  & &\frac{1}{2} &0&\frac{1}{2}\\
  & & &\frac{1}{2}&\frac{1}{2}
\end{pmatrix}.$$
From \cite[Theorem 8.1]{GLY24},
$$\Root\left(C^{+}_{2n}(x)\right)=\left\{\frac{1}{a^{k-1}}\cos\left(\frac{(2j-1)\pi}{k+2}\right):1\leq j\leq\frac{k+1}{2}\right\}.$$
Similarly,
$C^{-}_{2n}(x)$ is the characteristic polynomial of
$$\frac{1}{a^{k-1}}\begin{pmatrix}
 0 & \frac{1}{2} & \\
 \frac{1}{2} & 0 & \ddots\\
  & \ddots &\ddots&\frac{1}{2}\\
  & &\frac{1}{2} &0&\frac{1}{2}\\
  & & &\frac{1}{2}&-\frac{1}{2}
\end{pmatrix}.$$
From \cite[Theorem 8.1]{GLY24},
$$\Root\left(C^{-}_{2n}(x)\right)=\left\{\frac{1}{a^{k-1}}\cos\left(\frac{2j\pi}{k+2}\right):1\leq j\leq\frac{k+1}{2}\right\}.$$

\textbf{Case 2}, when $k=2n\geq 2$ is even. From Definition ~\ref{defn3.01}, we know that
$G^{+}_{2n+1}(x)$ is the characteristic polynomial of
$$\frac{1}{a^{k-1}}\begin{pmatrix}
 0 & \frac{1}{2} & \\
 \frac{1}{2} & 0 & \ddots\\
  & \ddots &\ddots&\frac{1}{2}\\
  & &\frac{1}{2} &0&\frac{\sqrt{2}}{2}\\
  & & &\frac{\sqrt{2}}{2}&0
\end{pmatrix}.$$
Knowing from \cite[Theorem 8.1]{GLY24},
$$\Root\left(G^{+}_{2n+1}(x)\right)=\left\{\frac{1}{a^{k-1}}\cos\left(\frac{(2j-1)\pi}{k+2}\right):1\leq j\leq\frac{k}{2}+1\right\}.$$

Similarly, $G^{-}_{2n+1}(x)$ is the characteristic polynomial of
$$\frac{1}{a^{k-1}}\begin{pmatrix}
 0 & \frac{1}{2} & \\
 \frac{1}{2} & 0 & \ddots\\
  & \ddots &\ddots&\frac{1}{2}\\
  & &\frac{1}{2} &0&\frac{1}{2}\\
  & & &\frac{1}{2}&0
\end{pmatrix}.$$
Knowing from \cite[Theorem 2.2]{KST99},
$$\Root\left(G^{-}_{2n+1}(x)\right)=\left\{\frac{1}{a^{k-1}}\cos\left(\frac{2j\pi}{k+2}\right):1\leq j\leq\frac{k}{2}\right\}.$$

Thus, by Theorem ~\ref{thm3.03},
$$\sigma_{p}(S_w\odot S_w^{*})=\left\{\left\{\frac{1}{a^{k-1}}\cos\left(\frac{(2j-1)\pi}{k+2}\right):k\ge 0,~1\leq j\leq \left \lfloor \frac{k+2}{2} \right \rfloor\right\}\right\},$$
$$\sigma_{p}(S_w\wedge S_w^{*})=\left\{\left\{\frac{1}{a^{k-1}}\cos\left(\frac{2j\pi}{k+2}\right):k\ge 1,~1\leq j\leq \left \lfloor \frac{k+1}{2} \right \rfloor\right\}\right\}.$$

When $k=4m-4,~j=m\ge1$, $\displaystyle \frac{1}{a^{k-1}}\cos\left(\frac{(2j-1)\pi}{k+2}\right)=0$, thus $\Dim\ker(S_w\odot S_w^{*})=\infty.$
Similarly, when $k=4m-2,~j=m\ge1$, $\displaystyle \frac{1}{a^{k-1}}\cos\left(\frac{2j\pi}{k+2}\right)=0$, thus $\Dim\ker(S_w\wedge S_w^{*})=\infty.$

Finally, when $a>1$, $S_w\odot S_w^{*}$ and $S_w\wedge S_w^{*}$ are compact self-adjoint operators. From Fredholm alternative theorem, we know
$$\sigma(S_w\odot S_w^{*})=\sigma_{p}(S_w\odot S_w^{*}),~\sigma(S_w\wedge S_w^{*})=\sigma_{p}(S_w\wedge S_w^{*}).$$

\end{pf}

\section{Symmetric tensor products of unilateral weighted shift operators and diagonal operators} \label{section4}

In this section, we will consider the symmetric tensor products of unilateral weighted shift operators and diagonal operators. This setting suggest us to work in the sequence space $\ell^{2}$ instead of $H^{2}(\mathbb{D})$. Let $\left\{e_{0},e_{1},\cdots\right\}$ be the standard basis of $\ell^{2}$, $\{\alpha_i\}_{i=0}^\infty$ be a bounded complex sequence and consider unilateral weighted shift operator $S_\alpha(e_{i})=\alpha_{i}e_{i+1},~i\in\mathbb{N}_0$. Clearly $S_\alpha^{*}(e_{i})=\alpha_{i-1}e_{i-1}$ for $i\ge 1~\textup{and}~S_\alpha^{*}(e_{0})=0.$

The following theorem is a generalization of \cite[Theorem 9.1]{GLY24}, which also rectifies some flaws of \cite[Theorem 9.1(b)]{GLY24}.
\begin{thm} \label{thm4.01}
Let $\{\alpha_i\}_{i=0}^\infty$ be a bounded complex sequence, $S_\alpha\in \mathcal{B}(\ell^2)$ be the unilateral weighted shift operator with $S_\alpha(e_{i})=\alpha_{i}e_{i+1}$, $\forall i\in \mathbb{N}_0$ and $M=\diag(\mu_{0},\mu_{1},\cdots)$ be a bounded diagonal operator on $\ell^{2}$. Assume that $\alpha_{i}\ne 0,~\forall i\in\mathbb{N}_0$.\\
\textup{(a)} If some $\mu_{i}=0$, then $0\in \sigma_{p}(S_\alpha\odot M)$, otherwise, $S_\alpha\odot M$ is injective.\\
\textup{(b)} $\sigma_{p}(S_\alpha\odot M)\subseteq\left\{0\right\}$.
\end{thm}

\begin{pf} (a) Note that for all $i,~j\ge 0$,
\begin{equation}
\begin{aligned}
(S_\alpha\odot M)(e_{i}\odot e_{j}) &=\frac{1}{2}(S_\alpha\otimes M+M\otimes S_\alpha)\left(\frac{1}{2}(e_{i}\otimes e_{j}+e_{j}\otimes e_{i})\right)\\
 &=\frac{1}{4}\left(S_\alpha(e_{i})\otimes M(e_{j})+M(e_{i})\otimes S_\alpha(e_{j})+S_\alpha(e_{j})\otimes M(e_{i})+M(e_{j})\otimes S_\alpha(e_{i})\right)\\
 &=\frac{1}{4}(\alpha_{i}\mu_{j}e_{i+1}\otimes e_{j}+\mu_{i}\alpha_{j}e_{i}\otimes e_{j+1}+\alpha_{j}\mu_{i}e_{j+1}\otimes e_{i}+\mu_{j}\alpha_{i}e_{j}\otimes e_{i+1})\\
 &=\frac{1}{2}\mu_{i}\alpha_{j}e_{i}\odot e_{j+1}+\frac{1}{2}\mu_{j}\alpha_{i}e_{i+1}\odot e_{j}.\nonumber
\end{aligned}
\end{equation}

If some $\mu_{i}=0$, then $0\in\sigma_{p}(S_\alpha\odot M)$ since
$$(S_\alpha\odot M)(e_{i}\odot e_{i})=0.$$

If $\mu_{i}\ne 0,~\textup{for}~\textup{all}~ i\ge 0$, let $\underset{0\leq i\leq j<\infty}{\sum}|b_{ij}|^{2}< \infty$ and let
$$v=2\sum_{0\leq i\leq j<\infty}b_{ij}e_{i}\odot e_{j}.$$
From Lemma ~\ref{lem2.08}, we know that $v\in \ell^{2}\odot\ell^{2}$.
\begin{equation}
\begin{aligned}
(S_\alpha\odot M)v &=(S_\alpha\odot M)\left(2\sum_{0\leq i\leq j<\infty}b_{ij}e_{i}\odot e_{j}\right)\\
 &=\sum_{0\leq i\leq j<\infty}b_{ij}(\mu_{j}\alpha_{i}e_{i+1}\odot e_{j}+\mu_{i}\alpha_{j}e_{i}\odot e_{j+1}).\nonumber
\end{aligned}
\end{equation}

For $(S_\alpha\odot M)v$, we consider the coefficient of $e_{k}\odot e_{l},~0\leq k\leq l<\infty$.
Actually, the coefficient is
$$
\left\{\begin{matrix}
 0 & \text{if}~~k=l=0, & \text{(1)}\\
 2\mu_{0}\alpha_{0}b_{0,0} & \text{if}~~k=0,~l=1, & \text{(2)}\\
 b_{0,l-1}\mu_{0}\alpha_{l-1} & \text{if}~~k=0,~l\ge 2, & \text{(3)}\\
 b_{k-1,k}\mu_{k}\alpha_{k-1} & \text{if}~~1\leq k=l, & \text{(4)}\\
 2b_{k,k}\mu_{k}\alpha_{k}+b_{k-1,k+1}\mu_{k+1}\alpha_{k-1} & \text{if}~~k\ge 1,~l=k+1, & \text{(5)}\\
 b_{k,l-1}\mu_{k}\alpha_{l-1}+b_{k-1,l}\mu_{l}\alpha_{k-1} & \text{if}~~k\ge 1,~l\ge k+2. &\text{(6)}
\end{matrix}\right.$$

We consider the upper triangular matrix of $\left\{b_{i,j}\right\}$,
$$
\begin{pmatrix}
b_{0,0}&b_{0,1}&b_{0,2}&b_{0,3}&b_{0,4}&b_{0,5}&\cdots\\
&b_{1,1}&b_{1,2}&b_{1,3}&b_{1,4}&b_{1,5}&\cdots\\
& &b_{2,2}&b_{2,3}&b_{2,4}&b_{2,5}&\cdots\\
& & &b_{3,3}&b_{3,4}&b_{3,5}&\cdots\\
& & & &b_{4,4}&b_{4,5}&\cdots\\
& & & & &b_{5,5}\\
& & & & & &\ddots
\end{pmatrix}.$$
If $(S_\alpha\odot M)v=0$, we will carry out inductive method to show that $v=0$.

Suppose that $P(k)$ is ``elements in row $k$ are all 0''. We will recursively use $\alpha_{n}, \mu_{n}\ne 0$, $\forall n\in \mathbb{N}_0$. From $(2), (3)$.
$$b_{0,0}=0,~b_{0,l-1}=0,~l\ge 2.$$
So the $P(0)$ is true. Now we assume that $P(n)$ is true, $0\leq n\leq N$.

$(4)$ tells us
$$b_{N+1,N+2}=0.$$
From $(5)$,
$$2b_{N+1,N+1}\mu_{N+1}\alpha_{N+1}+b_{N,N+2}\mu_{N+2}\alpha_{N}=0.$$
Since $P(N)$ is true, it follows that
\[b_{N+1,N+1}=0.\]
$(6)$ tells us
$$b_{N+1,l-1}\mu_{N+1}\alpha_{l-1}+b_{N,l}\mu_{l}\alpha_{N}=0,~l\ge N+3.$$
Again as $P(N)$ is true,
we have
\[b_{N+1,l}=0,~l\ge N+3.\] Consequently, $P(N+1)$ is true. Hence for all $ k,~P(k)$ is true.
Therefore $v=0$, and $S_\alpha\odot M$ is injective.\\

(b) Suppose that
$$\lambda \ne 0,~(S_\alpha\odot M)v=\lambda v,$$
where
$$v=2\sum_{0\leq i\leq j<\infty}b_{i,j}e_{i}\odot e_{j},\sum_{0\leq i\leq j<\infty}|b_{i,j}|^{2}<\infty.$$
Then we have
\begin{equation}
\begin{aligned}
0&=((S_\alpha\odot M)-\lambda I)v\\
&=2\sum_{0\leq i\leq j<\infty}b_{i,j}((S_\alpha\odot M)-\lambda I)(e_{i}\odot e_{j})\\
&=\sum_{0\leq i\leq j<\infty}b_{i,j}\mu_{j}\alpha_{i}e_{i+1}\odot e_{j}+\sum_{0\leq i\leq j<\infty}b_{i,j}\mu_{i}\alpha_{j}e_{i}\odot e_{j+1}-\sum_{0\leq i\leq j<\infty}\lambda b_{i,j}e_{i}\odot e_{j}.\nonumber
\end{aligned}
\end{equation}
Also, for $((S_\alpha\odot M)-\lambda I)v$, we consider the coefficient of $e_{k}\odot e_{l},~0\leq k\leq l<\infty$.
Actually, the coefficient is
$$
\left\{\begin{matrix}
 -\lambda b_{0,0} & \text{if}~~k=l=0, & \text{(1)}\\
 2\mu_{0}\alpha_{0}b_{0,0}-\lambda b_{0,1} & \text{if}~~k=0,~l=1, & \text{(2)}\\
 b_{0,l-1}\mu_{0}\alpha_{l-1}-\lambda b_{0,l} & \text{if}~~k=0,~l\ge 2, & \text{(3)}\\
 b_{k-1,k}\mu_{k}\alpha_{k-1}-\lambda b_{k,k} & \text{if}~~1\leq k=l, & \text{(4)}\\
 2b_{k,k}\mu_{k}\alpha_{k}+b_{k-1,k+1}\mu_{k+1}\alpha_{k-1}-\lambda b_{k,k+1} & \text{if}~~k\ge 1,~l=k+1, & \text{(5)}\\
 b_{k,l-1}\mu_{k}\alpha_{l-1}+b_{k-1,l}\mu_{l}\alpha_{k-1}-\lambda b_{k,l} & \text{if}~~k\ge 1,~l\ge k+2. &\text{(6)}
\end{matrix}\right.$$

We let $p(k)$ be
$$``b_{k,k+i}=0,~\textup{for}~\textup{all}~i\ge 0."$$
Since (1), (2), (3), we can know that
$$b_{0,0}=0,$$
$$2\mu_{0}\alpha_{0}b_{0,0}-\lambda b_{0,1}=0.$$
Recall that $\lambda\ne0$, hence,
$$b_{0,1}=0,~\textup{and}~b_{0,l-1}\mu_{0}\alpha_{l-1}-\lambda b_{0,l}=0,~l\ge 2.$$
Obviously, $p(0)$ is true.

Fix $N\in \mathbb{N}$, and we assume that $p(N-1)$ is true, then we have
$$b_{N-1,l}=0,~\textup{for}~\textup{all}~l\ge N-1.$$

(4) tells us $b_{N,N}=0$ and we can know from (5)
$$\lambda b_{N,N+1}=2b_{N,N}\mu_{N}\alpha_{N}+b_{N-1,N+1}\mu_{N+1}\alpha_{N-1}=0.$$
Then
$$b_{N,N+1}=0.$$

(6) says
$$\lambda b_{N,l}=b_{N,l-1}\mu_{N}\alpha_{l-1}+b_{N-1,l}\mu_{l}\alpha_{N-1}=b_{N,l-1}\mu_{N}\alpha_{l-1},~l\ge N+2.$$
Thus,
$$b_{N,l}=0,~\textup{for}~\textup{all}~l\ge N+2.$$
Hence $p(N)$ is true. Then we have $\textup{for}~\textup{all}~ k,~p(k)$ is true. Hence
$$v=0~,~\textup{thus}~\lambda \notin \sigma_{p}(S_\alpha\odot M).$$
Consequently, $\sigma_{p}(S_\alpha\odot M)\subseteq \left\{0\right\}.$

\end{pf}

We finally focus on $S_\alpha^{*}\odot M$ which is quite different from $S_\alpha\odot M$. Inspired by \cite[Problem 8]{GLY24}, the following theorem is a generalization of \cite[Theorem 9.2]{GLY24}.

\begin{thm} \label{thm4.02}
Let $\{\alpha_i\}_{i=0}^\infty$ be a bounded complex sequence, $S_\alpha\in \mathcal{B}(\ell^2)$ be the unilateral weighted shift operator with $S_\alpha(e_{i})=\alpha_{i}e_{i+1}$, $\forall i\in \mathbb{N}_0$ and $M=\diag(\mu_{0},\mu_{1},\cdots)$ be a bounded diagonal operator on $\ell^{2}$.  Then\\
\textup{(a)} $\displaystyle \frac{1}{\sqrt{2}}\underset{i\in \mathbb{N}}{\sup}\left(|\alpha_{i-1}||\mu_{i}|\right)\leq \|S_\alpha^{*}\odot M\|\leq \underset{i\in \mathbb{N}_0}{\sup}|\alpha_{i}|\underset{i\in \mathbb{N}_0}{\sup}|\mu_{i}|.$ Both inequalities are sharp.\\
\textup{(b)} $\displaystyle \left\{z:|z|<\frac{1}{2}|\mu_{0}|\inf_{j\in \mathbb{N}}|\alpha_{0}\cdots \alpha_{j-1}|^{\frac{1}{j}}\right\}\subset \sigma_{p}(S_\alpha^{*}\odot M)$, and $0\in \sigma_p(S_\alpha^*\odot M)$ if $\mu_0=0$.
\end{thm}

\begin{pf}
(a) First observe that
$$\|S_\alpha^{*}\|=\|S_\alpha\|=\underset{i\in \mathbb{N}_0}{\sup}|\alpha_{i}|.$$
From \cite[Theorem 3.4]{GLY24},
$$\|S_\alpha^{*}\odot M\|\leq \|S_\alpha^{*}\|\|M\|=\underset{i\in \mathbb{N}_0}{\sup}|\alpha_{i}|\|M\|=\underset{i\in \mathbb{N}_0}{\sup}|\alpha_{i}|\underset{i\in \mathbb{N}_0}{\sup}|\mu_{i}|.$$

We will show the left inequality (in (a)) in the following. Recall that $M(e_{i})=\mu_{i}e_{i}$ for $i\ge 0$. Then
$$(S_\alpha^{*}\odot M)(e_{i}\odot e_{j})=
\begin{cases}
 \frac{1}{2}(\mu_{j}\alpha_{i-1}e_{i-1}\odot e_{j}+\mu_{i}\alpha_{j-1}e_{i}\odot e_{j-1})&\text{if}~i,~j\ne 0,\\
 \frac{1}{2}\mu_{i}\alpha_{j-1}e_{i}\odot e_{j-1}&\text{if}~0=i<j,\\
 \frac{1}{2}\mu_{j}\alpha_{i-1}e_{i-1}\odot e_{j}&\text{if}~0=j<i,\\
 0&\text{if}~i=j=0.
\end{cases}$$
Thus
$$\|(S_\alpha^{*}\odot M)(e_{i}\odot e_{i})\|=\|\mu_{i}\alpha_{i-1}e_{i-1}\odot e_{i}\|=\frac{1}{\sqrt{2}}|\mu_{i}\alpha_{i-1}|.$$
Then, we have
$$\|S_\alpha^{*}\odot M\|\ge \sup_{i\in \mathbb{N}}\frac{1}{\sqrt{2}}|\mu_{i}\alpha_{i-1}|.$$

Let $M=I$, $S_\alpha(e_i)=e_{i+1}$, $\forall i\in \mathbb{N}_{0}$. Then $\|S_\alpha^{*}\|=r(S_\alpha^*)=1$. From \cite[Theorem 6.3]{GLY24},
$$\sigma(S_\alpha^{*}\odot I)=\frac{1}{2}(\sigma(S_\alpha^{*})+\sigma(S_\alpha^{*}))\subseteq \overline{\left\{z:|z|<\sup_{n\in \mathbb{N}_0}|\alpha_{n}|\right\}}=\overline{\mathbb{D}}.$$
Thus,
$$\|S_\alpha^{*}\odot I\|\ge 1.$$
Let $\mu_{i}=\delta_{i,1},~\alpha_{i}=\delta_{i,0},~i\ge 0$, where
$$\delta_{i,j}= \begin{cases}
0&~i\ne j,\\
1&~i=j.
\end{cases}$$
Then
$$\|S_\alpha^{*}\odot M\|=\|(S_\alpha^{*}\odot M)(e_{1}\odot e_{1})\|=\|\mu_{1}\alpha_{0}e_{0}\odot e_{1}\|=\frac{1}{\sqrt{2}}|\mu_{1}\alpha_{0}|=\sup_{i\in\mathbb{N}}\frac{1}{\sqrt{2}}|\mu_{i}\alpha_{i-1}|.$$
Hence, both inequalities are sharp.

(b) If $\mu_{0}=0$, then $(S_\alpha^{*}\odot M)(e_{1}\odot e_{0})=0$, thus $0\in \sigma_{p}(S_\alpha^{*}\odot M).$\\
Assume that
$$\mu_{0}\ne 0,~~\inf_{j\in \mathbb{N}}|\alpha_{0}\cdots \alpha_{j-1}|^{\frac{1}{j}}>0.$$
Fix $\lambda\in \mathbb{C}$ with
\[|\lambda|<\frac{1}{2}|\mu_{0}|\inf_{j\in \mathbb{N}}|\alpha_{0}\dots \alpha_{j-1}|^{\frac{1}{j}},\]
clearly, there exists $\beta\in (0,1)$ such that
\[\frac{|2\lambda|}{|\mu_0|\underset{j\in \mathbb{N}}{\inf}|\alpha_{0}\dots \alpha_{j-1}|^{\frac{1}{j}}}<\beta<1.\]
Thus,
$$\frac{|2\lambda|}{|\mu_{0}||\alpha_{0}\cdots \alpha_{j-1}|^{\frac{1}{j}}}<\beta,~~\forall j\in \mathbb{N}.$$
By Lemma ~\ref{lem2.08},
$$v\doteq e_0\odot e_0+\sum_{j=1}^{\infty}\frac{(2\lambda)^{j}}{\mu_{0}^{j}\alpha_{0}\cdots \alpha_{j-1}}e_{0}\odot e_{j}\in\ell^{2}\odot\ell^{2}.$$
Then since $(S_\alpha^*\odot M)(e_0\odot e_0)=0$,
\begin{equation}
\begin{aligned}
(S_\alpha^{*}\odot M)v &=\sum_{j=1}^{\infty}\frac{(2\lambda)^{j}}{\mu_{0}^{j}\alpha_{0}\cdots \alpha_{j-1}}(S_\alpha^{*}\odot M)(e_{0}\odot e_{j})\\
 &=\frac{1}{2}\frac{2\lambda}{\mu_0\alpha_0}\mu_0\alpha_0e_0\odot e_0+\frac{1}{2}\sum_{j=2}^{\infty}\frac{(2\lambda)^{j}}{\mu_{0}^{j}\alpha_{0}\cdots \alpha_{j-1}}\mu_{0}\alpha_{j-1}e_{0}\odot e_{j-1}\\
 &=\lambda e_0\odot e_0+\lambda\sum_{j=2}^{\infty}\frac{(2\lambda)^{j-1}}{\mu_{0}^{j-1}\alpha_{0}\cdots \alpha_{j-2}}e_{0}\odot e_{j-1}\\
 &=\lambda e_0\odot e_0+\lambda\sum_{n=1}^{\infty}\frac{(2\lambda)^{n}}{\mu_{0}^{n}\alpha_{0}\cdots \alpha_{n-1}}e_{0}\odot e_{n}\\
 &=\lambda v.\nonumber
\end{aligned}
\end{equation}
Hence, $\lambda \in \sigma_{p}(S_\alpha^{*}\odot M)$. The proof is complete.
%
%
\end{pf}

\begin{eg} \label{eg4.03}
Let $\alpha_{n}=1,~\textup{for}~\textup{all}~n\in \mathbb{N}_0,~M=I$, by Theorem ~\ref{thm4.02}, we have
$$\frac{1}{2}\mathbb{D}\subset\sigma_{p}(S_\alpha^{*}\odot I).$$
Since $\sigma_p(S_\alpha^*)=\mathbb{D}$,
for $\lambda\in\mathbb{D}$, one could choose a unit vector $e_\lambda$ such that $S_\alpha^{*}(e_{\lambda})=\lambda e_{\lambda}.$
Set
$f_{\lambda}=e_{\lambda}\otimes e_{\lambda}$. We have
$$(S_\alpha^{*}\odot I)f_{\lambda}=\lambda f_{\lambda},~\textup{for}~\textup{all}~\lambda\in\mathbb{D}.$$
Thus,
$$\frac{1}{2}\mathbb{D}\subset\mathbb{D}\subset\sigma_{p}(S_\alpha^{*}\odot I).$$
Consequently, $\frac{1}{2}\mathbb{D}\ne\sigma_{p}(S_\alpha^{*}\odot I).$ Thus, the inclusion relation in Theorem ~\ref{thm4.02} (b) is not sharp.
\end{eg}

The following Examples are the applications of Theorem ~\ref{thm4.02}.
%

\begin{eg} \label{eg4.05}(Dirichlet shift)
Let $\displaystyle \alpha_{n}=\sqrt{\frac{n+2}{n+1}},~ \textup{for}~\textup{all}~n\in \mathbb{N}_0$, by Theorem ~\ref{thm4.02}, we have
$$\left\{z:|z|<\frac{1}{2}|\mu_{0}|\right\}\cup\left\{0\right\}\subset \sigma_{p}(S_\alpha^{*}\odot M).$$
\end{eg}

When $\mu_{0}\ne 0$, from Theorem ~\ref{thm4.02} we know that
\[\underset{j\in \mathbb{N}}\inf|\alpha_{0}\dots \alpha_{j-1}|^{\frac{1}{j}}=\underset{j\in \mathbb{N}}\inf(j+1)^{\frac{1}{2j}}=1,\] thus \[\left\{z:|z|<\frac{1}{2}|\mu_{0}|\right\}\subset \sigma_{p}(S_\alpha^{*}\odot M).\]
 When $\mu_{0}=0$, then $(S_\alpha^{*}\odot M)(e_{1}\odot e_{0})=0$, thus $0\in \sigma_{p}(S_\alpha^{*}\odot M).$

\begin{eg} \label{eg4.06}(Bergman shift)
Let $\displaystyle \alpha_{n}=\sqrt{\frac{n+1}{n+2}},~\textup{for}~\textup{all}~n\in \mathbb{N}_0$, by Theorem ~\ref{thm4.02}, we have
$$\left\{z:|z|<\frac{\sqrt{2}}{4}|\mu_{0}|\right\}\cup\left\{0\right\}\subset \sigma_{p}(S_\alpha^{*}\odot M).$$
\end{eg}

When $\mu_{0}\ne 0$, from Theorem ~\ref{thm4.02} we know that
\[\underset{j\in \mathbb{N}}\inf|\alpha_{0}\dots \alpha_{j-1}|^{\frac{1}{j}}=\underset{j\in \mathbb{N}}\inf(\frac{1}{j+1})^{\frac{1}{2j}}=\frac{\sqrt{2}}{2},\]
thus \[\left\{z:|z|<\frac{\sqrt{2}}{4}|\mu_{0}|\right\}\subset \sigma_{p}(S_\alpha^{*}\odot M).\]
 When $\mu_{0}=0$, then $(S_\alpha^{*}\odot M)(e_{1}\odot e_{0})=0$, thus $0\in \sigma_{p}(S_\alpha^{*}\odot M).$

\section*{Acknowledgments}
The authors thank Professor Ji Youqing, Dou Rui, Luo Denghui, and Hu Zhaolong for many valuable discussions on the paper.

\bibliographystyle{plain}

\end{document}